\documentclass[a4paper,fleqn]{cas-sc}

\usepackage[authoryear,longnamesfirst]{natbib}
\usepackage{amsmath,amsfonts,amssymb}
\usepackage{enumitem}
\usepackage{microtype}
\usepackage{xcolor}

\newproof{pf}{Proof}
\renewcommand{\qed}{\hfill\ensuremath{\blacksquare}}
\newtheorem{Definition}{Definition}
\newtheorem{theorem}{Theorem}
\newtheorem{proposition}{Proposition}
\newtheorem{lemma}{Lemma}
\newtheorem{corollary}{Corollary}

\newtheorem{example}{Example}
\setlist[enumerate]{label=(\roman*),leftmargin=*,itemsep=2pt,topsep=4pt}

\newcommand{\PP}{\mathbb{P}}
\newcommand{\EE}{\mathbb{E}}
\newcommand{\R}{\mathbb{R}}
\newcommand{\st}{\mathrm{st}}
\newcommand{\calD}{\mathcal{D}}
\newcommand{\calE}{\mathcal{E}_A}

\hypersetup{
  pdftitle={Universal sampling rules for the signed friendship paradox},
  pdfauthor={Michail Louvaris and Samis Trevezas}
}

\begin{document}
\let\printorcid\relax

\shorttitle{Universal sampling rules for the signed friendship paradox}
\shortauthors{Louvaris and Trevezas}
\title[mode=title]{Universal sampling rules for the signed friendship paradox}

\author[1]{Michail Louvaris}
\cormark[1]
\author[2,3]{Samis Trevezas}
\address[1]{Department of Mathematics, KTH Royal Institute of Technology}
\address[2]{Department of Mathematics, National and Kapodistrian University of Athens,
Panepistimiopolis, Athens 15784, Greece}
\address[3]{MICS Laboratory, CentraleSup\'elec, Universit\'e Paris-Saclay,
3 Rue Joliot Curie, 91190 Gif-sur-Yvette, France}
\cortext[1]{Corresponding author. Email:
\href{mailto:louvaris@kth.se}{louvaris@kth.se}.}

\begin{abstract}
The friendship paradox describes a bias in mean degree induced by neighbour sampling. We study signed weighted networks, with transitions proportional to edge magnitudes and signed strength as the vertex score. We characterize all positive initial weights, depending on total strength and score, that guarantee a nonnegative one-step mean change on every connected graph. Without regularity assumptions, the same condition yields stochastic dominance at every fixed time. Continuous admissible rules are precisely total strength multiplied by a positive continuous nonincreasing function of the score. It is enough to test connected graphs with at most five vertices and two vertex types. The bound of five is optimal even for continuous rules; the proof uses a sharp four-vertex energy inequality. The complete score vector alone cannot determine a universally valid initial law with full support. Under marked local convergence and suitable integrability assumptions, we also obtain joint limits for local conditional mean score increments and compute an explicit limiting law and mean for independently signed random regular graphs. 
Under joint cut-norm convergence of the magnitude and signed kernels and uniformly bounded normalized strengths, we also obtain Wasserstein convergence of the sampling-weighted distributions of normalized local biases and convergence of their means.
\end{abstract}


\maketitle
\section{Introduction}

\citet{Feld1991} formulated the friendship paradox: friends have more
friends on average than individuals sampled uniformly from the population.
In the local formulation, one chooses an individual uniformly and then
one of their friends uniformly; on a finite undirected graph without
isolated vertices, the mean degree cannot decrease after this step.
Sampling uniformly from edge endpoints gives the degree-proportional
vertex law. Probabilistic, multistep and weighted versions are studied by
\citet{CaoRoss2016,Kramer2016,BerenhautJiang2019}.
\citet{Jackson2019} discusses consequences for perceptions and behaviour.

Negative interactions were considered by \citet{IyerCashore2016}
in their study of the downvoting paradox. \citet{GhasemianChristakis2023} introduced the
\emph{enmity paradox}, comparing the number of enemies of an individual
with those of their enemies, and studied mixed comparisons involving
friends and enemies. \citet{YangZhaoDerr2023} studied related comparisons
of positive and negative degrees. We consider signed weighted networks
in which neighbours are sampled in proportion to edge magnitudes and
the vertex score is signed strength. For an edge-weight matrix $A$,
these quantities are
\[
 d_i=\sum_j |A_{ij}|,
 \qquad
 s_i=\sum_j A_{ij}.
\]
Thus positive and negative weights can cancel in the score $s_i$
without reducing the total strength $d_i$ that determines encounter
probabilities.

Comparisons involving vertex attributes are studied by
\citet{EvtushenkoKleinberg2023}, and weighted versions by
\citet{EvtushenkoKleinberg2024}. Our question
concerns the initial sampling law: which positive weights $F(d_i,s_i)$
guarantee a nonnegative one-step mean change in signed score on every
finite connected signed weighted graph? For the walk considered here,
sampling proportionally to total strength is stationary and gives equality.
The corresponding observation for unsigned unweighted graphs is discussed by
\citet{Roga2026}.
We give a complete characterization of all universally admissible
rules, without regularity assumptions on $F$. The same condition
guarantees that the score at every time $k\ge1$ stochastically dominates
the initial score. Corollary~\ref{cor:regularity} shows that the
continuous admissible rules are precisely
\[
 F(a,r)=a g(r),
\]
where $g$ is positive, continuous and nonincreasing.

Universality can be tested on the family of connected signed weighted
graphs with at most five vertices and at most two vertex types
$(d_i,s_i)$. The bound of five is optimal even for continuous rules;
the sharpness proof uses a sharp energy inequality for four vertices.
We also show that the complete labelled score vector alone cannot
determine a universally valid initial law with full support. These
results identify both the admissible sampling rules and the information
needed to construct them.

For continuous admissible rules, Theorem~\ref{thm:marked-local} extends
the local-convergence results of
\citet{HazraEtAl2026a,HazraEtAl2026b} to signed weighted graphs and
nonuniform initial sampling. Under marked local convergence in
probability, uniform integrability of the sampling weights, and a
finite positive limiting normalizing constant, we obtain weak
convergence in probability of the sampling-weighted empirical
distributions of local conditional mean score increments, jointly
at finitely many fixed times. Additional uniform-integrability
assumptions yield convergence of their means. The limiting
distribution is obtained by reweighting the rooted local-limit
law by $d_o g(s_o)$ and then evaluating the corresponding local
bias vector.

Example~\ref{ex:signed-regular} applies this result to independently
signed random regular graphs. For the exponential sampling rule
$F_\beta(a,r)=a e^{-\beta r}$, $\beta\ge0$, we compute the one-step
limiting distribution and its mean explicitly. In particular, for
$d$-regular graphs with equally likely positive and negative edge
signs, the limiting mean is $(d-1)\tanh\beta$. 

Theorem \ref{thm:kernel-limit} provides a complementary limit result under kernel convergence, with applications to signed inhomogeneous independent-edge graphs in the dense regime and beyond. With the sampling rule applied to the normalized graph, joint cut-norm convergence of the magnitude and signed kernels and a uniform bound on normalized strengths holding with high probability yield joint convergence in probability in Wasserstein distance for the sampling-weighted local bias distributions at finitely many fixed times, together with convergence of their means.
\section{A universal classification}
A signed weighted graph on $V=\{1,\ldots,n\}$, $n\ge2$, is represented by a real symmetric
matrix $A$ with zero diagonal: $\{i,j\}$ is an edge exactly when $A_{ij}\ne0$. Assume that this
support graph is connected. Define the total strengths $d=(d_i)_{i\in V}$ and signed scores
$s=(s_i)_{i\in V}$ by
\begin{equation}\label{eq:model}
 d_i=\sum_{j\in V}|A_{ij}|,\quad s_i=\sum_{j\in V}A_{ij},\quad
 Z=\sum_{i\in V}d_i,\quad P_{ij}=\frac{|A_{ij}|}{d_i},\quad \pi_i=\frac{d_i}{Z}.
\end{equation}
Here $d_i>0$ and $|s_i|\le d_i$. Let $(X_k)_{k\ge0}$ be the Markov chain with transition matrix
$P$. For an initial law $\alpha$, write $\PP_\alpha$ and $\EE_\alpha$ for its probability and
expectation. Edge magnitudes govern encounters; signs enter the score $s(X_k)=s_{X_k}$.
Detailed balance, $\pi_iP_{ij}=\pi_jP_{ji}$, makes $\pi$ stationary.

Set $\calD=\{(a,r)\in\R^2:a>0,\ |r|\le a\}$. Fix one function $F:\calD\to(0,\infty)$, to be
used on every graph, and define
\begin{equation}\label{eq:initial}
 H(a,r)=\frac{F(a,r)}a,\qquad C_F=\sum_{j\in V}F(d_j,s_j),\qquad
 \alpha_i^F=\frac{F(d_i,s_i)}{C_F}.
\end{equation}
Write $\Delta_k(A,F)=\EE_{\alpha^F}[s(X_k)-s(X_0)]$, abbreviated to $\Delta_k$ when $A,F$ are fixed.
For a fixed finite graph, the initial law depends only on the finitely many
values $F(d_i,s_i)$, where $(d_i,s_i)$ determines the type of vertex $i$, so no measurability assumption on $F$ is needed here.  For real random variables $U,W$, the stochastic order $U\le_{\st}W$ means
$\PP(U>t)\le\PP(W>t)$ for every $t\in\R$.
Using $d_iP_{ij}=|A_{ij}|$ and summing the two directions of each edge gives
\begin{equation}\label{eq:edgeflow}
 C_F\Delta_1(A,F)=\sum_{i<j}|A_{ij}|
 \bigl(H(d_i,s_i)-H(d_j,s_j)\bigr)(s_j-s_i).
\end{equation}

\begin{theorem}\label{thm:classification}
For a fixed function $F:\calD\to(0,\infty)$, the following statements are equivalent.
\begin{enumerate}
\item $\Delta_1(A,F)\ge0$ for every finite connected signed weighted graph $A$.
\item $\Delta_1(A,F)\ge0$ for every connected signed weighted graph $A$ with at most five
  vertices and at most two vertex types.
\item For every $(a,u),(b,v)\in\calD$ with $u<v$, $H(a,u)\ge H(b,v)$.
\item For every finite connected signed weighted graph and every integer $k\ge1$, the chain
  started with law $\alpha^F$ satisfies $s(X_0)\le_{\st}s(X_k)$.
\end{enumerate}
The bound of five in (ii) is optimal, even among continuous $F$ when the four-vertex tests
allow all vertex types.
\end{theorem}

Condition (ii) concerns all real edge weights on the stated graph class.
Any rule that fails on some graph therefore has a counterexample with at most five vertices
and two types. We first prove the equivalence; the optimality of five is proved in Section~3.
\begin{pf}
We show (iii)$\Rightarrow$(iv)$\Rightarrow$(i)$\Rightarrow$(ii)$\Rightarrow$(iii).
For (iii)$\Rightarrow$(iv), fix a graph and an integer $k\ge1$. Write $P^k$ for the $k$-step
transition matrix and $Q^{(k)}_{ij}=\pi_i(P^k)_{ij}$. Detailed balance makes $Q^{(k)}$
symmetric. Write $h_i=\alpha_i^F/\pi_i$, so that $h_i=(Z/C_F)H(d_i,s_i)$. For any function
$\varphi$ on the finite score set,
\begin{equation}\label{eq:flow}
 \EE_{\alpha^F}\bigl[\varphi(s(X_k))-\varphi(s(X_0))\bigr]
 =\frac12\sum_{i,j}Q^{(k)}_{ij}(h_i-h_j)(\varphi(s_j)-\varphi(s_i)).
\end{equation}
Under (iii), each summand is nonnegative for nondecreasing $\varphi$. When $s_i=s_j$, the
summand vanishes regardless of $H(d_i,s_i)$ and $H(d_j,s_j)$. For each $t\in\R$, take
$\varphi(r)=\mathbf{1}_{\{r>t\}}$, the indicator of $(t,\infty)$, to obtain (iv). Taking
expectations in (iv) at $k=1$ gives (i), and restricting graphs gives (ii).

To prove (ii)$\Rightarrow$(iii), fix a type $(a,u)$ and put $p=(a+u)/2$, $q=(a-u)/2$. Choose
another type $(b,v)$ with
\begin{equation}\label{eq:budgets}
 p'=(b+v)/2>p/2,\qquad q'=(b-v)/2>q/2.
\end{equation}
Define a symmetric five-vertex matrix with zero diagonal by its upper-triangular entries:
\begin{equation}\label{eq:five}
\begin{aligned}
 A_{12}=A_{34}&=p'-p/2,& A_{23}=A_{15}=A_{45}&=p/2,\\
 A_{13}=A_{24}&=-(q'-q/2),& A_{14}=A_{25}=A_{35}&=-q/2.
\end{aligned}
\end{equation}
Vertices $1,\ldots,4$ have type $(b,v)$, and vertex $5$ has type $(a,u)$. By
\eqref{eq:budgets}, edges $12$, $24$, $43$ and $31$ form a cycle of nonzero weights.
Since $p+q=a>0$, vertex $5$ has an incident edge; the graph is connected, even when $p=0$
or $q=0$. Only edges incident to vertex $5$ join different scores. Their total magnitude
is $a$, so \eqref{eq:edgeflow} gives
\begin{equation}\label{eq:fivegap}
 \Delta_1=\frac{a(v-u)}{F(a,u)+4F(b,v)}\bigl(H(a,u)-H(b,v)\bigr).
\end{equation}
For arbitrary types $(a,u),(b,v)$ with $u<v$, choose $u<w<v$ and $c>\max(a,b)+|w|$. Then
$(c,w)$ satisfies \eqref{eq:budgets} relative to either type. Applying (ii) with vertex $5$
of type $(a,u)$ and the other vertices of type $(c,w)$ gives $H(a,u)\ge H(c,w)$.
Applying it with vertex $5$ of type $(b,v)$ gives $H(c,w)\ge H(b,v)$, since $w-v<0$.
The two inequalities prove
(iii), including boundary types $|u|=a$ or $|v|=b$.
\qed
\end{pf}

To describe possible discontinuities, write $g(r+)=\lim_{t\downarrow r}g(t)$ and
$g(r-)=\lim_{t\uparrow r}g(t)$ for the one-sided limits of a nonincreasing function $g$.
\begin{corollary}\label{cor:regularity}
The equivalent conditions of Theorem~\ref{thm:classification} hold if and only if there is a
positive nonincreasing $g:\R\to(0,\infty)$ such that
\begin{equation}\label{eq:limits}
 a\,g(r+)\le F(a,r)\le a\,g(r-),\qquad (a,r)\in\calD.
\end{equation}
Consequently, strength dependence beyond the factor $a$ can occur only at countably many score
levels. If $F$ is continuous, the conditions are equivalent to $F(a,r)=a\,g(r)$ for a positive
continuous nonincreasing $g$.
\end{corollary}
\begin{pf}
For necessity, assume Theorem~\ref{thm:classification}(iii). Then $g(r)=H(|r|+1,r)$ is positive
and nonincreasing. For $(a,r)\in\calD$, condition (iii) gives $H(a,r)\ge g(t)$ if $t>r$ and
$H(a,r)\le g(t)$ if $t<r$. Taking one-sided limits gives \eqref{eq:limits}.
Conversely, for $u<v$, monotonicity gives $g(u+)\ge g(v-)$, so \eqref{eq:limits} implies (iii).
A monotone function has at most countably many discontinuities, and at each continuity point
\eqref{eq:limits} forces $H(a,r)=g(r)$. If $F$ is continuous on $\calD$, the chosen $g$ is
continuous, so \eqref{eq:limits} gives $F(a,r)=ag(r)$ on $\calD$.
Conversely, $F(a,r)=ag(r)$ with positive continuous nonincreasing $g$ satisfies
(iii).
\qed
\end{pf}

Theorem~\ref{thm:classification}(iv) compares each fixed time with time zero; it does not
assert monotonicity between successive times. Reversibility and (iii) give, for every
$k\ge1$ and $u<v$,
\begin{equation}\label{eq:orientation}
 \PP_{\alpha^F}(s(X_0)=u,s(X_k)=v)
 \ge \PP_{\alpha^F}(s(X_0)=v,s(X_k)=u).
\end{equation}
Put $D_k=s(X_k)-s(X_0)$. For each $t>0$, summing \eqref{eq:orientation} over score pairs
with difference $t$ gives $\PP_{\alpha^F}(D_k=t)\ge\PP_{\alpha^F}(D_k=-t)$.
Comparing the tails of these finite-support laws gives $D_k\ge_{\st}-D_k$. Under $F(a,r)=ag(r)$, the ratio of the two probabilities in
\eqref{eq:orientation}, when the denominator is positive, is $g(u)/g(v)$. For this factored
rule, \eqref{eq:flow} at $k=1$ gives $\Delta_1=0$ if and only if $g(s_i)$ is constant on
the connected graph. If $g$ is strictly decreasing, this is equivalent to all scores
coinciding.

\section{Why four vertices are insufficient}
For any real symmetric zero-diagonal matrix $A$, connected or not, define $d$ and $s$ as in
\eqref{eq:model}. For $x,y\in\R^n$, write
\begin{equation*}
 \calE(x,y)=\sum_{i<j}|A_{ij}|(x_i-x_j)(y_i-y_j),\qquad
 \calE(x)=\calE(x,x).
\end{equation*}
The following inequality holds for graphs with at most four vertices and is used to establish
the optimality of the bound in Theorem~\ref{thm:classification}.
\begin{lemma}\label{lem:energy}
For every real symmetric zero-diagonal matrix of order at most four,
\begin{equation}\label{eq:energy}
 \calE(d)\le3\calE(s).
\end{equation}
The constant $3$ is sharp. No finite constant can replace $3$ for all matrices of order five.
\end{lemma}
\begin{pf}
If the matrix has order less than four, adjoin zero rows and columns. Write
$a=A_{12}$, $a'=A_{34}$, $b=A_{13}$, $b'=A_{24}$, $c=A_{14}$ and $c'=A_{23}$.
Set $x=(a-a')/2$, $y=(b-b')/2$ and $z=(c-c')/2$. By the reverse triangle inequality,
there are $\lambda,\mu,\nu\in[-1,1]$ satisfying
$(|a|-|a'|)/2=\lambda x$, $(|b|-|b'|)/2=\mu y$ and $(|c|-|c'|)/2=\nu z$;
take $\lambda=0$ when $x=0$, and define $\mu,\nu$ in the same way. Set
$p=|a|+|a'|$, $q=|b|+|b'|$, $r=|c|+|c'|$ and
\begin{equation*}
 L=(q+r)x^2+(p+r)y^2+(p+q)z^2.
\end{equation*}
Since $p\ge2|x|$, $q\ge2|y|$ and $r\ge2|z|$, the arithmetic--geometric mean inequality gives
\begin{equation*}
L\ge 2\bigl[x^2(|y|+|z|)+y^2(|x|+|z|)+z^2(|x|+|y|)\bigr]\ge12|xyz|.
\end{equation*}
Direct calculation gives
\begin{align*}
 \calE(s)/4&=L+4(\lambda+\mu+\nu)xyz,\\
 \calE(d)/4&=(q+r)\lambda^2x^2+(p+r)\mu^2y^2+(p+q)\nu^2z^2
                  +12\lambda\mu\nu xyz.
\end{align*}
The function $\lambda+\mu+\nu-\lambda\mu\nu$ is affine in each coordinate and has
absolute value at most $2$ at every corner of $[-1,1]^3$. Therefore
\begin{equation*}
 (3\calE(s)-\calE(d))/4
 \ge 2L+12(\lambda+\mu+\nu-\lambda\mu\nu)xyz
 \ge2L-24|xyz|\ge0.
\end{equation*}
The star with $A_{12}=A_{13}=1$, $A_{14}=-1$ and all other edges absent has $d=(3,1,1,1)$,
$s=(1,1,1,-1)$, and energies $12$ and $4$, so the constant $3$ is attained. For the five-vertex
graph in \eqref{eq:five} with $(a,u)=(1,0)$ and $(b,v)=(4,0)$, all scores vanish while
$\calE(d)=9$. Hence no finite constant works for five vertices.
\qed
\end{pf}

\begin{proposition}\label{prop:sharp}
There is a positive continuous $F$ for which $\Delta_1(A,F)\ge0$ on every connected signed
weighted graph $A$ with at most four vertices, but $\Delta_1(A,F)<0$ on a graph with five
vertices.
\end{proposition}
\begin{pf}
We start with a nonincreasing function of the score and add a small term supported near
strength $3/2$ and score $0$. For $t\in\R$, let $t_+=\max(t,0)$ and define
\begin{equation}\label{eq:bump}
 \begin{aligned}
 c(r)&=\max(-1,\min(r,1)),\qquad
 f(a,r)=(1-|2a-3|)_+(1-4|r|)_+,\\
 F_*(a,r)&=a\bigl(2-c(r)+10^{-3}f(a,r)\bigr).
 \end{aligned}
\end{equation}
This function is continuous and strictly positive, since $2-c(r)\ge1$. On a graph with at
most four vertices, put
$f_i=f(d_i,s_i)$, write $c(s)=(c(s_i))_{i\in V}$ and $\mathbf{f}=(f_i)_{i\in V}$, and set
$B=\calE(c(s),s)\ge0$, $T=\calE(\mathbf{f},s)$. By \eqref{eq:edgeflow},
$C_{F_*}\Delta_1(A,F_*)=B-10^{-3}T$. We show that
\begin{equation}\label{eq:perturbation}
 |T|\le684B.
\end{equation}
If all $f_i$ vanish, then $T=0$. Otherwise choose $i$ with $f_i>0$; it satisfies $1<d_i<2$ and
$|s_i|<1/4$. The total magnitudes of its positive and negative incident edges are $(d_i+s_i)/2$
and $(d_i-s_i)/2$, both greater than $3/8$. Each sign class contains at most two edges, so an
edge from each class has magnitude greater than $3/16$.

If $B<3/256$, the scores of these two neighbours lie within $1/4$ of $s_i$. Indeed, any
difference of at least $1/4$ contributes at least $(3/16)(1/4)^2$ to $B$: since
$|s_i|<1/4$, clipping the other score outside $[-1,1]$ still leaves a clipped difference
of at least $1/4$. Three scores therefore belong to $(-1/2,1/2)$. If a fourth score
satisfies $|s_j|\ge1$, every edge incident to that vertex contributes at least $|A_{j\ell}|/4$
to $B$, implying $B\ge d_j/4\ge1/4$. This is impossible. Thus all scores lie in $(-1,1)$ and
$B=\calE(s)$. The bound
\begin{equation*}
 |f(a,r)-f(b,t)|\le2|a-b|+4|r-t|
\end{equation*}
and Cauchy--Schwarz, followed by Lemma~\ref{lem:energy}, give
\begin{equation*}
 |T|\le2\sqrt{\calE(d)\calE(s)}+4\calE(s)
 \le(2\sqrt3+4)B.
\end{equation*}

If $B\ge3/256$, consider first the edges contributing to $T$ with $|s_i-s_j|\le1$. The
sum of their absolute contributions is at most
\begin{equation*}
 \sum_{\substack{i<j\,:\, f_i>0\text{ or }f_j>0}}|A_{ij}|
 \le\sum_{i:f_i>0}d_i\le8,
\end{equation*}
since $0\le f_i\le1$, there are at most four vertices with $f_i>0$, and each has $d_i<2$.
On each remaining edge that contributes to $T$, an endpoint with $f_i>0$ has score in
$(-1/4,1/4)$, so $|c(s_i)-c(s_j)|\ge3/4$. Its absolute contribution to $T$ is at most $4/3$
times its contribution to $B$. Hence
\begin{equation*}
 |T|\le8+\tfrac43B\le\left(\tfrac{2048}3+\tfrac43\right)B=684B,
\end{equation*}
which proves \eqref{eq:perturbation}. Since $684\cdot10^{-3}<1$, this estimate gives
$C_{F_*}\Delta_1(A,F_*)\ge(316/1000)B\ge0$.

Finally, take the five-vertex graph \eqref{eq:five} where vertex $5$ has type
$(a,u)=(1,-1/2000)$ and vertices $1,\ldots,4$ have type $(b,v)=(3/2,0)$. Condition
\eqref{eq:budgets} holds. Writing $H_*(a,r)=F_*(a,r)/a$, we have
\begin{equation*}
 H_*(1,-1/2000)=\frac{4001}{2000}
 <\frac{2001}{1000}=H_*(3/2,0).
\end{equation*}
Formula \eqref{eq:fivegap} gives $\Delta_1=-1/56026000<0$.
\qed
\end{pf}

Proposition~\ref{prop:sharp} also proves the optimality assertion of
Theorem~\ref{thm:classification}.

\section{What information must the initial rule use?}
A sampling law has full support if it gives positive probability to every vertex. We now allow
it to depend arbitrarily on the complete labelled score vector $s$, rather than through a
vertexwise function $F$.
\begin{proposition}\label{prop:information}
No full-support initial law determined solely by the labelled score vector guarantees a
nonnegative one-step mean gap on every connected signed weighted graph. For one fixed score
vector, support and edge-sign pattern, every law assigning positive mass to both score groups
gives both signs of the gap as the edge magnitudes vary.
\end{proposition}
\begin{pf}
Fix $u<v$, the score vector $(u,u,u,u,v,v,v,v)$, and $\varepsilon>0$. Choose
$a>|u|$ and $b>|v|$ so that $(a+u)/2>\varepsilon$ and $(b+v)/2>\varepsilon$. Set
$p=(a+u)/2$, $q=(a-u)/2$, $p'=(b+v)/2$, $q'=(b-v)/2$, and let $I_4$ be the
$4\times4$ identity matrix. Define
\begin{equation*}
 M(p,q)=\begin{pmatrix}0&p&-q&0\\p&0&0&-q\\-q&0&0&p\\0&-q&p&0\end{pmatrix},\qquad
 A_\varepsilon=\begin{pmatrix}M(p-\varepsilon,q)&\varepsilon I_4\\
 \varepsilon I_4&M(p'-\varepsilon,q')\end{pmatrix}.
\end{equation*}
The two groups have respective types $(a,u)$ and $(b,v)$, and the graph is connected. For a
fixed initial law $\alpha$ with mass $\rho\in(0,1)$ on the first group, its one-step gap is
\begin{equation*}
 \EE_\alpha[s(X_1)-s(X_0)]=\varepsilon(v-u)\left(\frac\rho a-\frac{1-\rho}b\right).
\end{equation*}
Keeping $\varepsilon$ fixed, increasing $a$ with $b$ fixed makes this negative; increasing $b$
with $a$ fixed makes it positive. The scores, support and signs remain unchanged. A rule based
only on that score vector must assign the same initial law in both cases.
\qed
\end{pf}

Without the full-support restriction, choosing only minimum-score vertices always gives a nonnegative gap.

\section{Local limits and an application}\label{sec:local}

For unweighted graphs, \citet[Theorem~2.2]{HazraEtAl2026a} establish
convergence of the empirical one-step friendship-bias distribution under
local convergence in probability. The corresponding fixed-time result in
the locally tree-like setting appears in
\citet[Theorem~2(a)]{HazraEtAl2026b}. We give a marked version that
incorporates the sampling weights characterized above and yields joint
convergence at finitely many fixed times, followed by an explicit
application.

\subsection*{Marked local convergence}

Let $\mathcal G_*^{\mathrm m}$ be the space of rooted isomorphism classes
of connected, locally finite simple graphs with symmetric, nonzero real
edge marks. Write $H=(G,A,o)$ for an element of this space. The marked
ball of radius $r$ is the induced subgraph on vertices at graph distance
at most $r$ from $o$, with its root and edge marks. A sequence converges
in the marked local topology if, for each fixed $r$, these balls
eventually have the same support up to a root-preserving isomorphism,
and the corresponding marks converge uniformly over that finite ball.
This is the usual rooted-network topology; see
\citet[Section~2]{AldousLyons2007}. If a finite graph is disconnected,
a rooted graph denotes the component containing its root.

We use the following notion of marked local convergence in probability.
\begin{Definition}\label{def:local}
Let $(G_n,A_n)$ be random signed weighted graphs on $[n]$, and let $U_n$
be conditionally uniform on $[n]$ given $(G_n,A_n)$. Put
\[
 L_n=\frac1n\sum_{i=1}^n\delta_{(G_n,A_n,i)}.
\]
We say that $(G_n,A_n)$ converges locally in probability to a
deterministic probability law $\mu$ on $\mathcal G_*^{\mathrm m}$ if,
for every bounded continuous
$h:\mathcal G_*^{\mathrm m}\to\mathbb R$,
\[
 L_n(h)
 =
 \mathbb E[h(G_n,A_n,U_n)\mid G_n,A_n]
 \xrightarrow{\mathbb P}\mu(h).
\]
\end{Definition}

This is the marked analogue of
\citet[Definition~2.1(b)]{HazraEtAl2026a}. We use uniform integrability
in its usual sense: $(Y_n)_n$ is uniformly integrable if
\[
 \lim_{M\to\infty}\sup_n
 \mathbb E\bigl[
   |Y_n|\mathbf1_{\{|Y_n|>M\}}
 \bigr]=0.
\]

\subsection*{Local biases and sampling weights}

Fix a continuous function $g:\mathbb R\to(0,\infty)$ and put
$F(a,r)=a g(r)$. In this section graphs may be disconnected and may
have isolated vertices; we extend the rule by $F(0,0)=0$. For
$H=(G,A,o)$, set
\[
 d_v=\sum_{u\sim v}|A_{vu}|,\qquad
 s_v=\sum_{u\sim v}A_{vu},\qquad
 q_v=d_vg(s_v),\qquad
 q(H)=q_o.
\]
Let $P_{vu}=|A_{vu}|/d_v$ when $d_v>0$, and make isolated vertices
absorbing. For $k\ge1$, define
\[
 b_k(H)
 =
 (P^k s)_o-s_o
 =
 \mathbb E_o^H[s(X_k)-s_o],
 \qquad
 \boldsymbol b_K(H)
 =
 (b_1(H),\ldots,b_K(H)).
\]
The expectation is over the walk on the fixed marked graph $H$.
Thus $b_k$ is a conditional mean score increment. At a nonisolated
root,
\[
 b_1
 =
 \frac1{d_o}\sum_{u\sim o}|A_{ou}|s_u-s_o;
\]
at an isolated root, $q=b_k=0$.

On $(G_n,A_n)$, write $q_{i,n}$, $b_{k,i,n}$ and
$\boldsymbol b_{K,i,n}$ for the corresponding quantities at vertex
$i$. Define
\[
 C_{F,n}=\sum_{i=1}^n q_{i,n}.
\]
When $C_{F,n}>0$, put
\[
 \alpha^F_{i,n}
 =
 \frac{q_{i,n}}{C_{F,n}},
 \qquad
 \Delta_{k,n}
 =
 \Delta_k(A_n,F)
 =
 \frac1{C_{F,n}}
 \sum_{i=1}^n q_{i,n}b_{k,i,n}.
\]
For a rooted marked law $\mu$ with $c=\mu(q)\in(0,\infty)$, define
the sampling-biased law
\[
 \widehat\mu_g(dH)
 =
 \frac{q(H)}c\,\mu(dH).
\]
Here $T_\#\eta$ denotes the pushforward of $\eta$ under $T$. This
change of root agrees with the finite-graph sampling rule:
\[
 \frac{qL_n}{L_n(q)}
 =
 \sum_{i=1}^n
 \alpha^F_{i,n}\delta_{(G_n,A_n,i)}
 \qquad\text{on }\{C_{F,n}>0\}.
\]

The following theorem gives joint limits at finitely many fixed times, under both uniform and weighted root sampling.
\begin{theorem}\label{thm:marked-local}
Suppose $(G_n,A_n)$ converges locally in probability in the marked
sense of Definition~\ref{def:local} to a deterministic law $\mu$ on
$\mathcal G_*^{\mathrm m}$. For every fixed integer $K\ge1$,
\begin{equation}\label{eq:unweighted-joint}
 \nu_{n,K}
 :=
 \frac1n\sum_{i=1}^n\delta_{\boldsymbol b_{K,i,n}}
 \Longrightarrow
 \nu_K
 :=
 (\boldsymbol b_K)_\#\mu
 \quad\text{weakly in probability on }\mathbb R^K.
\end{equation}
If $(q_{U_n,n})_n$ is uniformly integrable and
$c:=\mu(q)\in(0,\infty)$, then also
\begin{equation}\label{eq:weighted-joint}
 \boldsymbol\mu^g_{n,K}
 :=
 \frac1{C_{F,n}}
 \sum_{i=1}^n
 q_{i,n}\delta_{\boldsymbol b_{K,i,n}}
 \Longrightarrow
 \boldsymbol\mu^g_K
 :=
 (\boldsymbol b_K)_\#\widehat\mu_g
 \quad\text{weakly in probability on }\mathbb R^K.
\end{equation}
If, in addition,
$(q_{U_n,n}|b_{k,U_n,n}|)_n$ is uniformly integrable for each
$1\le k\le K$, then $\mu(q|b_k|)<\infty$ and
\begin{equation}\label{eq:joint-mean}
 (\Delta_{1,n},\ldots,\Delta_{K,n})
 \xrightarrow{\mathbb P}
 \frac1c
 \bigl(\mu(qb_1),\ldots,\mu(qb_K)\bigr).
\end{equation}
Every limiting mean in \eqref{eq:joint-mean} is nonnegative when
$g$ is nonincreasing. Coordinate projection in
\eqref{eq:weighted-joint} gives the marginal limit
$(b_k)_\#\widehat\mu_g$ for each $k\le K$, with mean
$\mu(qb_k)/c$ under the additional integrability assumption.

All uniform-integrability assumptions include the randomness of
the graph, its marks and the conditionally uniform root $U_n$.
On $\{C_{F,n}=0\}$, set
$\boldsymbol\mu^g_{n,K}=\delta_{\boldsymbol0}$ and
$\Delta_{k,n}=0$.
\end{theorem}

\begin{pf}
We first check that the local biases are continuous.
For a nonisolated root, expansion over length-$k$ walks gives
\[
 b_k(H)
 =
 \sum_{\substack{
       v_0=o\\
       v_\ell\sim v_{\ell+1},\ 0\le\ell<k
      }}
 \left(
   \prod_{\ell=0}^{k-1}
   \frac{|A_{v_\ell v_{\ell+1}}|}{d_{v_\ell}}
 \right)
 (s_{v_k}-s_o).
\]
Repeated vertices and edges are allowed. Local finiteness makes
the sum finite, and all its terms are determined by the marked
radius-$(k+1)$ ball. Under marked local convergence these finite
supports eventually agree, their marks converge, and every
denominator converges to a positive strength. At an isolated
limiting root, approximating roots are eventually isolated.
Hence each $b_k$ is continuous, as is $q$ by continuity of $g$.
In particular, $\boldsymbol b_K$ is a continuous radius-$(K+1)$
observable.

We next extend local convergence to a continuous real-valued
function $\Psi$ for which
$\Psi_n:=\Psi(G_n,A_n,U_n)$ is uniformly integrable. Define
\[
 T_M(x)=\max\{-M,\min\{x,M\}\}.
\]
For each fixed $M$,
\[
 L_n(T_M\circ\Psi)
 \xrightarrow{\mathbb P}
 \mu(T_M\circ\Psi).
\]
Boundedness and local convergence applied to $|\Psi|\wedge M$
give
\[
 \mu(|\Psi|\wedge M)
 =
 \lim_{n\to\infty}
 \mathbb E[|\Psi_n|\wedge M]
 \le
 \sup_n\mathbb E|\Psi_n|
 <
 \infty.
\]
Monotone convergence implies $\mu(|\Psi|)<\infty$. For
$\varepsilon>0$, conditional uniformity of $U_n$ and Markov's
inequality yield
\[
 \mathbb P\bigl(
   |L_n(\Psi)-L_n(T_M\circ\Psi)|>\varepsilon
 \bigr)
 \le
 \varepsilon^{-1}
 \mathbb E\bigl[
   |\Psi_n|\mathbf1_{\{|\Psi_n|>M\}}
 \bigr].
\]
The right-hand side tends to zero uniformly in $n$ as
$M\to\infty$; also
$\mu(|\Psi-T_M\circ\Psi|)\to0$. Consequently,
\begin{equation}\label{eq:unbounded-local}
 L_n(\Psi)
 \xrightarrow{\mathbb P}
 \mu(\Psi).
\end{equation}

For the empirical laws, take $f\in C_b(\mathbb R^K)$. Continuity gives
\[
 L_n(f\circ\boldsymbol b_K)
 \xrightarrow{\mathbb P}
 \mu(f\circ\boldsymbol b_K).
\]
Under the first uniform-integrability assumption,
\eqref{eq:unbounded-local} also applies to $q$ and to
$q(f\circ\boldsymbol b_K)$, since the latter is bounded in
absolute value by $\|f\|_\infty q$. Thus
\[
 \frac{C_{F,n}}n
 =
 L_n(q)
 \xrightarrow{\mathbb P}
 c>0,
 \qquad
 \frac{L_n(q(f\circ\boldsymbol b_K))}{L_n(q)}
 \xrightarrow{\mathbb P}
 \frac{\mu(q(f\circ\boldsymbol b_K))}{c}.
\]
In particular, $\mathbb P(C_{F,n}=0)\to0$. A countable
convergence-determining class in $C_b(\mathbb R^K)$ proves
\eqref{eq:unweighted-joint} and \eqref{eq:weighted-joint}.
Under the remaining assumptions, apply
\eqref{eq:unbounded-local} to $q|b_k|$ and $qb_k$ and divide by
$L_n(q)$. This gives finite limiting absolute moments and
convergence in probability of each mean, hence
\eqref{eq:joint-mean} because $K$ is fixed.

Finally, suppose that $g$ is nonincreasing. Reversibility gives
\[
 d_{i,n}(P_n^k)_{ij}
 =
 d_{j,n}(P_n^k)_{ji},
\]
including vertices in different components and isolated vertices
under the stated convention. Symmetrisation therefore yields,
when $C_{F,n}>0$,
\[
 C_{F,n}\Delta_{k,n}
 =
 \frac12
 \sum_{i,j}
 d_{i,n}(P_n^k)_{ij}
 \bigl(g(s_{i,n})-g(s_{j,n})\bigr)
 (s_{j,n}-s_{i,n})
 \ge0
\]
for nonincreasing $g$. The convention on $\{C_{F,n}=0\}$ and
\eqref{eq:joint-mean} complete the proof.
\qed
\end{pf}

\paragraph{One-step and classical special cases.}
For $K=1$, write
$\mu_n^g=\boldsymbol\mu^g_{n,1}$ and
$\mu_g=(b_1)_\#\widehat\mu_g$. The theorem gives both the
uniform-root limit $(b_1)_\#\mu$ and the sampling-weighted
limit $\mu_g$, where
\[
 \mu_g(B)
 =
 \frac{\mu(q\mathbf1_{\{b_1\in B\}})}{\mu(q)},
 \qquad B\in\mathcal B(\mathbb R).
\]
The marked radius-two law determines this one-step limit.
On an unsigned unweighted graph, take
$g(r)=1/\max\{1,r\}$. Then $s_v=d_v$ and
$q_v=\mathbf1_{\{d_v>0\}}$, so sampling is uniform over
nonisolated vertices. If $p_0=\mu(d_o=0)<1$, the classical
uniform-root bias law is
\[
 (b_1)_\#\mu
 =
 p_0\delta_0+(1-p_0)\mu_g.
\]
Thus it agrees with $\mu_g$ when $p_0=0$; compare
\citet{HazraEtAl2026a}.

For general fixed $K$, the limiting coordinates refer to the
same root and need not be independent. The means compare each
fixed time with time zero. The theorem concerns fixed $K$ and
does not assert convergence at depths growing with $n$ or
monotonicity between successive times.

For independently signed random regular graphs, the one-step limiting law can be computed explicitly.
\begin{example}\label{ex:signed-regular}
Fix an integer $d\ge3$, $p\in[0,1]$ and $\beta\ge0$. Let $G_n$
be uniform among simple $d$-regular graphs on $[n]$, along
integers $n>d$ with $nd$ even. Independently of $G_n$, assign
independent signs to its edges, positive with probability $p$,
and use these signs as the edge weights. Choose
\[
 g_\beta(s)=e^{-\beta s},
 \qquad
 F_\beta(a,r)=a e^{-\beta r}.
\]
Then $d_{i,n}=d$ and
\[
 q_{i,n}=d e^{-\beta s_{i,n}},
 \qquad
 \alpha^{F_\beta}_{i,n}
 =
 \frac{e^{-\beta s_{i,n}}}
      {\sum_j e^{-\beta s_{j,n}}}.
\]

The marked local limit is the rooted infinite $d$-regular tree
with independent edge signs of the same law. For the underlying
unmarked graphs, this follows from
\citet[Theorem~2.17]{VanDerHofstad2024}. To verify the marking
step, fix a radius $r$ and a bounded function $h$ of the marked
radius-$r$ ball. Let $R_{n,r}$ be the fraction of vertices whose
unmarked radius-$r$ ball is not the corresponding tree ball.
Then $R_{n,r}\to0$ in probability, and
\[
 \left|
   \mathbb E[L_n(h)\mid G_n]-\mu(h)
 \right|
 \le
 2\|h\|_\infty R_{n,r}.
\]
Conditionally on $G_n$, the values of $h$ at roots at distance
greater than $2r$ depend on disjoint sets of edge signs and are
independent. For fixed $d,r$, there are only $O(n)$ ordered
pairs of roots within distance $2r$, so
\[
 \operatorname{Var}(L_n(h)\mid G_n)=O(n^{-1}).
\]
Conditional Chebyshev's inequality gives
$L_n(h)\to\mu(h)$ in probability. Since the degree bound is
fixed and the mark space $\{-1,+1\}$ is finite, bounded
continuous functions can be uniformly approximated on this
compact rooted space by finite-radius functions. This proves
marked local convergence in the sense of
Definition~\ref{def:local}.

All integrability assumptions in
Theorem~\ref{thm:marked-local} hold automatically, since
\[
 d e^{-\beta d}\le q_{i,n}\le d e^{\beta d},
 \qquad
 |b_{k,i,n}|\le2d
 \quad(k\ge1).
\]
The theorem therefore gives joint convergence and convergence
of means at every finite collection of fixed times. The
one-step limit has an explicit form. Define
\[
 M_\beta
 =
 p e^{-\beta}+(1-p)e^\beta,
 \qquad
 p_\beta
 =
 \frac{p e^{-\beta}}{M_\beta},
 \qquad
 c=\mu(q)=d M_\beta^d,
\]
and let
\[
 B_\beta\sim\operatorname{Bin}(d,p_\beta),
 \qquad
 C\sim\operatorname{Bin}(d(d-1),p)
\]
be independent. Then
\begin{equation}\label{eq:regular-limit-law}
 \mu_{g_\beta}
 =
 \mathcal L\!\left(
   \frac2d\bigl(C-(d-1)B_\beta\bigr)
 \right),
 \qquad
 \Delta_1(A_n,F_\beta)
 \xrightarrow{\mathbb P}
 2(d-1)(p-p_\beta).
\end{equation}

Indeed, on the limiting tree let $S$ be the sum of the $d$
signs incident to the root and let $T$ be the sum of the
$d(d-1)$ signs from its neighbours to the next generation.
Under the uniform-root law, $S$ and $T$ are independent, and
\[
 b_1=\frac{T-(d-1)S}{d}.
\]
The sampling tilt $d e^{-\beta S}$ changes the root-edge signs
to independent signs with positive probability $p_\beta$ and
leaves $T$ independent with its original law. Thus
$S=2B_\beta-d$ and $T=2C-d(d-1)$ under
$\widehat\mu_{g_\beta}$, which proves
\eqref{eq:regular-limit-law}. The limiting mean is nonnegative
because $p_\beta\le p$, and is strictly positive when
$0<p<1$ and $\beta>0$. For symmetric signs, $p=1/2$, it becomes
\[
 \Delta_1(A_n,F_\beta)
 \xrightarrow{\mathbb P}
 (d-1)\tanh\beta.
\]
\end{example}
\section{A law of large numbers under kernel convergence}
\label{sec:kernel-limit}

Kernel convergence provides a counterpart to the sparse local-limit
result. Related dense limits for an unsigned triangle bias are obtained
by \citet[Section~1.3.3]{BhattacharyaEtAl2026}. Here both the
sampling-weighted bias distribution and its mean converge at every
fixed walk length.

Let $A_n$ be a symmetric signed weight matrix on $[n]$, with zero
diagonal, and let $(\lambda_n)$ be a deterministic sequence of positive numbers. Write
\[
 d_i=\sum_j|A_{ij}|,
 \qquad
 s_i=\sum_jA_{ij},
 \qquad
 (P_n)_{ij}=\frac{|A_{ij}|}{d_i}
 \quad\text{when }d_i>0,
\]
with isolated vertices absorbing. Apply the fixed rule
$F(a,r)=ag(r)$ to the normalized graph
$\widetilde A_n=A_n/\lambda_n$, where
$g:\mathbb R\to(0,\infty)$ is continuous. The initial probabilities
and normalized local biases are
\begin{equation}\label{eq:kernel-normalization}
 \widetilde\alpha_{i,n}
 =\frac{d_i g(s_i/\lambda_n)}
        {\sum_jd_jg(s_j/\lambda_n)},
 \qquad
 \widetilde b_{k,i,n}
 =\frac{(P_n^ks)_i-s_i}{\lambda_n}.
\end{equation}
Thus isolated vertices receive zero sampling weight. Define
\[
 \widetilde\Delta_{k,n}
 =\sum_i\widetilde\alpha_{i,n}\widetilde b_{k,i,n},
 \qquad
 \widetilde\mu^g_{n,K}
 =\sum_i\widetilde\alpha_{i,n}
   \delta_{(\widetilde b_{1,i,n},\ldots,
            \widetilde b_{K,i,n})}.
\]
If the normalizing sum is zero, set these quantities to $0$ and
$\delta_{\boldsymbol 0}$, respectively.

For $I_i=((i-1)/n,i/n]$, define the magnitude and signed kernels
\[
 U_n(x,y)=\frac{n}{\lambda_n}|A_{ij}|,
 \qquad
 W_n(x,y)=\frac{n}{\lambda_n}A_{ij},
 \qquad (x,y)\in I_i\times I_j.
\]
For an integrable kernel $H$, write
\[
 \|H\|_\square
 =\sup_{S,T\subseteq[0,1]}
   \left|\int_{S\times T}H(x,y)\,dx\,dy\right|,
 \qquad
 (Hf)(x)=\int_0^1H(x,y)f(y)\,dy,
\]
where the supremum is over measurable sets.

The next result gives convergence of the weighted bias distributions
under a common relabelling of the two kernels.
\begin{theorem}\label{thm:kernel-limit}
Suppose that, after a common measure-preserving relabelling of
$U_n,W_n$,
\begin{equation}\label{eq:kernel-assumptions}
 \|U_n-U\|_\square+\|W_n-W\|_\square
 \xrightarrow{\mathbb P}0,
 \qquad
 \mathbb P\!\left(
   \max_i\frac{d_i}{\lambda_n}\le C
 \right)\longrightarrow1
\end{equation}
for some $C<\infty$ and deterministic symmetric integrable kernels
$U,W$ satisfying $U\ge|W|$ and $\iint U>0$.

Put $d=U\mathbf1$, $s=W\mathbf1$, and define
\[
 (\mathcal Pf)(x)
 =
 \begin{cases}
  (Uf)(x)/d(x), & d(x)>0,\\
  f(x), & d(x)=0.
 \end{cases}
\]
Let
\[
 b_k=\mathcal P^ks-s,
 \qquad
 \boldsymbol b_K=(b_1,\ldots,b_K),
 \qquad
 c_g=\int_0^1d(x)g(s(x))\,dx.
\]
Then $0<c_g<\infty$ and, for every fixed $K\ge1$,
\begin{equation}\label{eq:kernel-law}
 W_1\!\left(
   \widetilde\mu^g_{n,K},
   (\boldsymbol b_K)_\#
   \left[\frac{d(x)g(s(x))}{c_g}\,dx\right]
 \right)
 \xrightarrow{\mathbb P}0.
\end{equation}
Here $T_\#\nu$ denotes the pushforward of $\nu$ under $T$, and
$W_1$ is the Wasserstein distance of order one on $\mathbb R^K$
for the $\ell^1$ distance.
In particular, for every fixed $k\ge1$,
\begin{equation}\label{eq:kernel-mean}
 \widetilde\Delta_{k,n}
 \xrightarrow{\mathbb P}
 \frac{
   \int_0^1d(x)g(s(x))
       \bigl((\mathcal P^ks)(x)-s(x)\bigr)\,dx
 }{c_g}.
\end{equation}
Every limiting mean is nonnegative if $g$ is nonincreasing.
\end{theorem}

\begin{pf}
Set $d_n=U_n\mathbf1$, $s_n=W_n\mathbf1$, and let $\mathcal P_n$
be the corresponding kernel walk. The estimate
\[
 \|Hf\|_1\le4\|H\|_\square\|f\|_\infty
\]
gives
$\|d_n-d\|_1+\|s_n-s\|_1\to0$ in probability.
On $\{\|d_n\|_\infty\le C\}$,
$\int(d-C)_+\le\|d-d_n\|_1$, hence $d\le C$ almost everywhere.
We work on this event, whose probability tends to one.
Then $|s_n|,|s|\le C$, and their walk averages are bounded by $C$.

Symmetry gives the weighted contraction
\[
 \int d_n|\mathcal P_nh|
 \le \int d_n|h|.
\]
Moreover,
\[
 d_n(\mathcal P_n-\mathcal P)f
 =(U_n-U)f-(d_n-d)\mathcal Pf,
\]
including zero-degree points. Consequently,
\begin{equation*}
 \int d_n|\mathcal P_nf_n-\mathcal Pf| \
 \le{} \
 \int d_n|f_n-f|
\  + \ 4\|U_n-U\|_\square\|f\|_\infty 
 \ + \ \|d_n-d\|_1\|\mathcal Pf\|_\infty.
\end{equation*}
Starting with
$\int d_n|s_n-s|\le C\|s_n-s\|_1$ and iterating this inequality
shows that, for every fixed $K$,
\[
 \int d_n
 \|\boldsymbol b_{n,K}-\boldsymbol b_K\|_1
 \xrightarrow{\mathbb P}0,
 \qquad
 \boldsymbol b_{n,K}
 =(\mathcal P_n^ks_n-s_n)_{k=1}^K.
\]
Put $G=\max_{[-C,C]}g$. Uniform continuity of $g$ on $[-C,C]$ gives
$\|g(s_n)-g(s)\|_1\to0$ in probability on the event under consideration, and
\begin{equation*}
 \|d_ng(s_n)-dg(s)\|_1
 \ \le \ G\|d_n-d\|_1+C\|g(s_n)-g(s)\|_1
 \ \xrightarrow{\mathbb P}\ 0.
\end{equation*}

Since $g$ has a positive minimum on $[-C,C]$ and $\int d>0$,
we have $0<c_g<\infty$.

Write $q_n=d_ng(s_n)$, $q=dg(s)$ and $c_n=\int q_n$.
Then $c_n\to c_g>0$ in probability, so $\PP(c_n>0)\to1$.
On this event, $\|q_n/c_n-q/c_g\|_1\to0$ in probability, and the finite weighted law equals
$(\boldsymbol b_{n,K})_\#(q_n\,dx/c_n)$.
Coupling the two bias vectors at the same point, and then comparing the root densities, gives
\begin{equation*}
 W_1\!\left(\widetilde\mu^g_{n,K},
        (\boldsymbol b_K)_\#(q\,dx/c_g)\right)
 \ \le \ \
 \frac{G}{c_n}\int d_n
      \|\boldsymbol b_{n,K}-\boldsymbol b_K\|_1\,dx
 +2CK\left\|\frac{q_n}{c_n}-\frac{q}{c_g}\right\|_1
 \ \xrightarrow{\mathbb P}\ 0.
\end{equation*}
The second term uses the diameter $4CK$ of $[-2C,2C]^K$ and
$\|\rho-\eta\|_{\mathrm{TV}}=\tfrac12\|\rho-\eta\|_1$ for probability densities.
The coordinate maps are $1$-Lipschitz, so \eqref{eq:kernel-mean} follows.
For nonincreasing $g$, each finite-graph mean is nonnegative by reversibility,
and the same is true of its limit.\qed
\end{pf}

\paragraph{Application: signed inhomogeneous independent-edge graphs.}
Suppose the $A_{ij}$, $i<j$, are independent and satisfy
the following bounds, with constants $M,C_0>0$ independent of $n$:
\[
 |A_{ij}|\le M,
 \qquad
 \max_i\sum_j\mathbb E|A_{ij}|\le C_0\lambda_n,
 \qquad
 \liminf_{n\to\infty}\frac{\lambda_n}{\log n}>0.
\]
If the expected kernels $\mathbb EU_n,\mathbb EW_n$ converge
jointly in cut norm to $U,W$ as in the theorem, then
\eqref{eq:kernel-assumptions} holds.
By Bernstein's inequality \citep[Section~2.8]{BoucheronLugosiMassart2013},
for $L>C_0$,
\begin{equation*}
 \PP\!\left(\max_i d_i>L\lambda_n\right)
 \le n\exp\{-\kappa_L\lambda_n\},\qquad
 \kappa_L=\frac{(L-C_0)^2}{2M\{C_0+(L-C_0)/3\}}.
\end{equation*}
Since $\kappa_L\to\infty$ as $L\to\infty$, the lower bound on
$\lambda_n/\log n$ gives the required strength bound for sufficiently large $L$.
For a fixed vertex cut, the centred sum is a sum of independent unordered-edge
contributions, with uniformly bounded magnitudes and variance $O(n\lambda_n)$.
A union bound over at most $4^n$ vertex cuts therefore gives, for each fixed $\varepsilon>0$,
\begin{equation*}
 \PP\!\left(\|U_n-\EE U_n\|_\square>\varepsilon\right)
 \le 2\cdot4^n\exp\{-c_\varepsilon n\lambda_n\}
 \longrightarrow0
\end{equation*}
for some $c_\varepsilon>0$. The same argument applies to $W_n$.
Together with the assumed convergence of the expected kernels, these bounds prove
\eqref{eq:kernel-assumptions}.
Thus the theorem applies to bounded signed marks with
endpoint-dependent distributions, allowing positive, negative,
or zero signed means.

For an explicit example, let $n$ be even, take two vertex classes of size $n/2$,
and retain each edge independently with probability $\rho_n$.
Present edges are positive within the first class, negative within
the second, and have independent fair signs between classes.
For $\lambda_n=n\rho_n$ with
$\liminf_n\lambda_n/\log n>0$, the limiting magnitude kernel is
$U=1$, while the signed kernel has block values
\[
 W=
 \begin{pmatrix}
  1&0\\
  0&-1
 \end{pmatrix}.
\]
Hence $d=1$, $s=\pm\tfrac12$, and $\mathcal P^ks=0$ for every
$k\ge1$. The limiting marginal bias law is
\[
 \frac{
  g(1/2)\delta_{-1/2}+g(-1/2)\delta_{1/2}
 }{
  g(1/2)+g(-1/2)
 },
\]
and, for every fixed $k\ge1$,
\[
 \widetilde\Delta_{k,n}
 \xrightarrow{\mathbb P}
 \frac{
  g(-1/2)-g(1/2)
 }{
  2\{g(1/2)+g(-1/2)\}
 }.
\]
In particular, for $g(r)=e^{-\beta r}$ with $\beta\ge0$, the
limiting mean is $\tfrac12\tanh(\beta/2)$.

The normalization in \eqref{eq:kernel-normalization} is part of the
model: the fixed sampling rule is applied to $A_n/\lambda_n$.
For fixed raw-score weights $d_i g(s_i)$, additional control of
$g$ at growing scores is needed.
The two limiting kernels must also be kept distinct:
one may have $|W|<U$ because signs cancel in the score without
cancelling encounter magnitudes.

\paragraph{Discussion.}
The rules $F(a,r)=ag(r)$ have a sampling interpretation. Sample an edge proportionally to
$|A_{ij}|$ and select an endpoint uniformly, obtaining law $\pi$. For $M\ge\max_i g(s_i)$,
retain endpoint $i$ with probability $g(s_i)/M$; the retained law is $\alpha^F$.
The classification concerns a rule used on every graph; a fixed graph may admit additional
initial laws. The necessity proof uses arbitrary real edge weights. Restricting the weights
to $\{-1,0,1\}$ requires a separate characterization on the admissible integer values of
total and signed strength.

\paragraph{Acknowledgements.}
We thank Frank den Hollander for his May 2025 talk on the friendship paradox at the Department
of Mathematics of the National and Kapodistrian University of Athens and for a subsequent
discussion with ST, which together inspired this letter.

%\paragraph{Declaration of generative AI and AI-assisted technologies in the manuscript preparation process.}
The authors used ChatGPT (OpenAI) for a revision of an initial draft of the article that resulted in improvements in the organization and presentation of the final manuscript. The authors reviewed
and revised the material and take full responsibility for the results, proofs,
references and final text of this article.

\bibliographystyle{cas-model2-names}
\bibliography{signed_sampling_references}

@article{Feld1991,
  author = {Feld, Scott L.},
  title = {Why Your Friends Have More Friends Than You Do},
  journal = {American Journal of Sociology},
  year = {1991}, volume = {96}, number = {6}, pages = {1464--1477},
  doi = {10.1086/229693}
}

@article{CaoRoss2016,
  author = {Cao, Y. and Ross, Sheldon M.},
  title = {The friendship paradox},
  journal = {The Mathematical Scientist},
  year = {2016}, volume = {41}, number = {1}, pages = {61--64}
}

@article{Kramer2016,
  author = {Brown Kramer, Josh and Cutler, Jonathan and Radcliffe, A. J.},
  title = {The multistep friendship paradox},
  journal = {The American Mathematical Monthly},
  year = {2016}, volume = {123}, number = {9}, pages = {900--908},
  doi = {10.4169/amer.math.monthly.123.9.900}
}

@article{BerenhautJiang2019,
  author = {Berenhaut, Kenneth S. and Jiang, Hongyi},
  title = {The friendship paradox for weighted and directed networks},
  journal = {Probability in the Engineering and Informational Sciences},
  year = {2019}, volume = {33}, number = {1}, pages = {136--145},
  doi = {10.1017/S0269964818000050}
}

@article{Jackson2019,
  author = {Jackson, Matthew O.},
  title = {The friendship paradox and systematic biases in perceptions and social norms},
  journal = {Journal of Political Economy},
  year = {2019}, volume = {127}, number = {2}, pages = {777--818},
  doi = {10.1086/701031}
}

@article{EvtushenkoKleinberg2023,
  author = {Evtushenko, Anna and Kleinberg, Jon},
  title = {Node-based generalized friendship paradox fails},
  journal = {Scientific Reports},
  year = {2023}, volume = {13}, pages = {2074},
  doi = {10.1038/s41598-023-29268-7}
}

@article{EvtushenkoKleinberg2024,
  author = {Evtushenko, Anna and Kleinberg, Jon},
  title = {A comprehensive generalization of the {Friendship Paradox} to weights and attributes},
  journal = {Scientific Reports},
  year = {2024}, volume = {14}, pages = {13730},
  doi = {10.1038/s41598-024-63167-9}
}

@article{GhasemianChristakis2023,
  author = {Ghasemian, Amir and Christakis, Nicholas A.},
  title = {The enmity paradox},
  journal = {Scientific Reports},
  year = {2023}, volume = {13}, pages = {20040},
  doi = {10.1038/s41598-023-47167-9}
}

@inproceedings{IyerCashore2016,
  author = {Iyer, Shankar and Cashore, Massey},
  title = {Friendship Paradoxes and the {Quora} Downvoting Paradox},
  booktitle = {JSM Proceedings, Section on Statistics in Marketing},
  year = {2016}, pages = {52--67},
  publisher = {American Statistical Association},
  address = {Alexandria, VA},
  url = {https://ww2.amstat.org/meetings/proceedings/2016/data/assets/pdf/389489.pdf}
}

@misc{Roga2026,
  author = {Roga, Wojciech},
  title = {Friendship paradox disappears under degree biased network sampling},
  year = {2026},
  note = {arXiv:2603.16337},
  url = {https://arxiv.org/abs/2603.16337}
}

@misc{YangZhaoDerr2023,
  author = {Yang, Catherine and Zhao, Yuying and Derr, Tyler},
  title = {The Friendship Paradox: An Analysis on Signed Social Networks with Positive and Negative Links},
  year = {2023},
  howpublished = {KDD 2023 Undergraduate Consortium},
  note = {Nonarchival consortium paper},
  url = {https://www.kdd.org/kdd2023/ph-d-and-undergraduate-consortiums/index.html}
}

@article{AldousLyons2007,
  author  = {Aldous, David and Lyons, Russell},
  title   = {Processes on Unimodular Random Networks},
  journal = {Electronic Journal of Probability},
  year    = {2007},
  volume  = {12},
  pages   = {1454--1508},
  doi     = {10.1214/EJP.v12-463}
}

@article{HazraEtAl2026a,
  author  = {Hazra, Rajat Subhra and den Hollander, Frank
             and Parvaneh, Azadeh},
  title   = {The friendship paradox for sparse random graphs},
  journal = {Probability Theory and Related Fields},
  year    = {2026},
  volume  = {195},
  pages   = {109--131},
  doi     = {10.1007/s00440-025-01365-w}
}

@article{HazraEtAl2026b,
  author  = {Hazra, Rajat Subhra and den Hollander, Frank
             and Parvaneh, Azadeh},
  title   = {The multi-level friendship paradox for sparse random graphs},
  journal = {Stochastic Processes and their Applications},
  year    = {2026},
  volume  = {195},
  pages   = {104873},
  doi     = {10.1016/j.spa.2026.104873}
}

@book{VanDerHofstad2024,
  author    = {van der Hofstad, Remco},
  title     = {Random Graphs and Complex Networks. Volume 2},
  series    = {Cambridge Series in Statistical and Probabilistic Mathematics},
  number    = {54},
  publisher = {Cambridge University Press},
  year      = {2024},
  doi       = {10.1017/9781316795552}
}

@article{BhattacharyaEtAl2026,
  author  = {Bhattacharya, Bishakh and Gadhiwala, Nitya
             and {den Hollander}, Frank and Jain, Pradeeptha R
             and Subramanya, Tejashree},
  title   = {The Triangle Friendship Paradox},
  journal = {Sankhya A},
  year    = {2026},
  doi     = {10.1007/s13171-026-00458-5},
  note    = {Published online 18 September 2026}
}

@book{BoucheronLugosiMassart2013,
  author    = {Boucheron, St{\'e}phane and Lugosi, G{\'a}bor
               and Massart, Pascal},
  title     = {Concentration Inequalities:
               A Nonasymptotic Theory of Independence},
  publisher = {Oxford University Press},
  address   = {Oxford},
  year      = {2013},
  doi       = {10.1093/acprof:oso/9780199535255.001.0001}
}

\end{document}